\documentclass[11pt]{amsart}

\usepackage{amsxtra,amssymb,amsmath,amscd,url,listings,stmaryrd}
\usepackage[alphabetic,initials]{amsrefs}
\usepackage{scalerel,stackengine}
\stackMath
\usepackage[utf8]{inputenc}
\usepackage{eucal}
\usepackage{fullpage}
\usepackage{scrtime}
\usepackage[colorlinks]{hyperref}
\usepackage{xcolor}
\usepackage{datetime}
\usepackage{graphicx}

\shortdate

\renewcommand{\leq}{\leqslant}
\renewcommand{\geq}{\geqslant}
\newcommand\widecheck[1]{%
\savestack{\tmpbox}{\stretchto{%
  \scaleto{%
    \scalerel*[\widthof{\ensuremath{#1}}]{\kern-.4pt\bigwedge\kern-.4pt}%
    {\rule[-\textheight/2]{1ex}{\textheight}}
  }{\textheight}%
}{0.5ex}}%
\stackon[2pt]{#1}{\scalebox{-1}{\tmpbox}}%
}

\newcommand{\expo}{1/260}

\numberwithin{equation}{section}

\def\stacksum#1#2{{\stackrel{{\scriptstyle #1}}
{{\scriptstyle #2}}}}

\newcommand{\sym}{\mathrm{sym}}

\newcommand{\Cc}{\mathbf{C}}

\newcommand{\Aa}{\mathbf{A}}

\newcommand{\Zz}{\mathbf{Z}}

\newcommand{\Rr}{\mathbf{R}}

\newcommand{\Qq}{\mathbf{Q}}

\newcommand{\Fq}{{\mathbf{F}_q}}

\newcommand{\Fqt}{{\mathbf{F}^\times_q}}

\newcommand{\KL}{\mathcal{K}\ell}

\newcommand{\mods}[1]{\,(\mathrm{mod}\,{#1})}

\newcommand{\what}{\widehat}

\newcommand{\ra}{\rightarrow}

\DeclareMathOperator{\Kl}{\mathrm{Kl}}

\newcommand{\eps}{\varepsilon}
\renewcommand{\rho}{\varrho}

\DeclareMathOperator{\GL}{GL}

\DeclareMathSymbol{\gena}{\mathord}{letters}{"3C}
\DeclareMathSymbol{\genb}{\mathord}{letters}{"3E}

\def\intc{\frac{1}{2i\pi}\mathop{\int}\limits}

\theoremstyle{plain}
\newtheorem{theorem}{Theorem}[section]
\newtheorem*{theorem*}{Theorem}
\newtheorem{lemma}[theorem]{Lemma}

\theoremstyle{remark}

\theoremstyle{definition}

\newtheorem{remark}[theorem]{Remark}

\newcommand{\mcF}{\mathcal{F}}

\newcommand{\mcK}{\mathcal{K}}

\newcommand{\vphi}{\varphi}

\renewcommand{\geq}{\geqslant}
\renewcommand{\leq}{\leqslant}
\renewcommand{\Re}{\mathfrak{Re}\,}

\newcommand{\ov}[1]{\overline{#1}}

\newcommand\sumsum{\mathop{\sum\sum}\limits}

\begin{document}
\title{Rankin--Selberg coefficients in large arithmetic progressions}
\author{Emmanuel Kowalski}\address{D-MATH, ETH Z\"urich, R\"amistrasse 101, CH-8092 Z\"urich, Switzerland}
\email{kowalski@math.ethz.ch}

\author{Yongxiao Lin}
\address{Data Science Institute, Shandong University, Jinan 250100, China}
\email{yongxiao.lin@sdu.edu.cn}

\author{Philippe Michel}
\address{EPFL/MATH/TAN, Station 8, CH-1015 Lausanne, Switzerland }
\email{philippe.michel@epfl.ch}

\begin{abstract} 
  Let $(\lambda_f(n))_{n\geq 1}$ be the Hecke eigenvalues of either a holomorphic Hecke eigencuspform or a Hecke--Maass cusp form $f$. We prove that, for any fixed $\eta>0$, under the Ramanujan--Petersson conjecture for $\rm GL_2$ Maass forms the Rankin--Selberg coefficients $(\lambda_f(n)^2)_{n\geq 1}$ admit a level of distribution $\theta=2/5+\expo-\eta$ in arithmetic progressions.
\end{abstract}

  \thanks{Ph.\ M. was partially supported by the SNF (grant 200021\_197045).  \today\ \currenttime}

\maketitle

\hfill{\em In memory of Jingrun Chen}
\section{Introduction}
In a series of papers \cites{FKM1,KLMS,LMS} we studied the absence of correlations between the coefficients of certain automorphic $L$-functions and trace functions of prime moduli. 

More precisely, given $q$ a prime number, let
$$K:\Fq=\Zz/q\Zz\to \Cc$$ be the trace function associated to a suitable $\ell$-adic middle extension sheaf $\mcF$  on the affine line $\Aa^1_{\Fq}$, geometrically irreducible and pure of weight $0$; this implies in particular that the supnorm satisfies
$$\|K\|_\infty\leq C(\mcF)$$
where $C(\mcF)$ denotes the analytic conductor of $\mcF$, a numerical invariant attached to the Galois representation underlying $\mcF$. We now view $K$ as a $q$-periodic function on $\Zz$ via the obvious projection.

Let
$$L(\pi,s)=\sum_{n\geq 1}\frac{\lambda_\pi(n)}{n^s}=\prod_p L(\pi_p,s),\ \Re s>1$$
be an automorphic $L$-function of some degree $d\geq 2$ (normalized so that $\Re s=1/2$ is the critical line). For $V$ a smooth, compactly supported function on $\Rr_{>0}$, we consider the problem of obtaining non-trivial bounds for the correlation sums
\begin{equation}\label{nontrivialgeneric}
S_V(K;X)=\sum_{n\geq 1}\lambda_\pi(n)K(n)V\left(\frac{n}X\right)\ll X^{1-\eta}\hbox{  as $q,X\ra\infty$;}	
\end{equation}
here $\eta>0$ is some positive constant and the above bound depends implicitly on $\pi$, $V$ and $C(\mcF)$. Under relatively mild conditions on $\mcF$ it is not too difficult to obtain non-trivial bounds like \eqref{nontrivialgeneric} as long as
	$$X\geq q^{d/2+\delta}$$
for some $\delta>0$ (with the exponent $\delta$ depending on $\eta$) and so the first challenging range is 
\begin{equation}\label{convexrange}
X\asymp q^{d/2}.
\end{equation}
This range is called the {\em convexity range} as it corresponds to the critical range for the subconvexity problem in the large $q$-aspect for the twisted $L$-function $L(\pi\times\chi,s)$ for $\chi\mods q$ a non-trivial Dirichlet character.

In the three papers mentioned above, a non-trivial bound \eqref{nontrivialgeneric} was obtained for $X$ at and below the convexity range. Specifically
\begin{itemize}
\item \cite{FKM1} considered the situation where $L(\pi,s)$ is the standard $L$-function of a $\GL_{2,\Qq}$ automorphic representation (the $L$-function of a Hecke eigenform) and obtained (under some suitable assumptions on $\mcF$) \eqref{nontrivialgeneric} as long as
\begin{equation}\label{shortestd2}X\geq q^{1-1/4+\delta},\ \delta>0.	
\end{equation}

\item \cite{KLMS} considered the situation where $L(\pi,s)$ is the standard $L$-function of a $\GL_{3,\Qq}$ automorphic representation (of level $1$) and obtained (again under some mild assumptions on $\mcF$)  \eqref{nontrivialgeneric} as long as
\begin{equation}\label{shortestd3}
	X\geq q^{3/2-1/6+\delta},\ \delta>0.
\end{equation}

\item \cite{LMS} considered the situation where $L(\pi,s)$ is the  Rankin--Selberg $L$-function attached to a pair $(\vphi,f)$ of
$\GL_{3,\Qq}$ and $\GL_{2,\Qq}$ automorphic forms (both of level $1$). More precisely  $\lambda_\pi$ is given by
\begin{equation}
	\label{RS23def}\lambda_\pi(n)=\sum_{mr^2=n}\lambda_\vphi(m,r)\lambda_f(m).
\end{equation} In that case \eqref{nontrivialgeneric} can be obtained for $K$ a trace function  associated with a suitably ``good" sheaf $\mcK$ (see \cite{LMS}*{\S 1} for the definition of the goodness) as long as
\begin{equation}\label{shortestd6}X\geq q^{3-1/4+\delta},\ \delta>0.	
\end{equation}
\end{itemize}

\subsection{Applications to large arithmetic progressions}

There is some genuine interest in trying to obtain \eqref{nontrivialgeneric} for even shorter ranges. One possible motivation is the study of the distribution of the sequences $(\lambda_\pi(n))_{n\leq X}$ in an arithmetic progression $n\equiv a\mods q$ for $(a,q)=1$ when $q\leq X$: set
\begin{equation}
    \label{Deltadef}
    \Delta(\lambda_\pi,X,a;q):=\mathop{\sum_{n\geq 1}}_{n\equiv a\bmod q}\lambda_\pi(n)V\left(\frac{n}{X}\right)-\frac{1}{\vphi(q)}
     \sum_\stacksum{n\geq 1}{(n,q)=1}\lambda_\pi(n)V\left(\frac{n}{X}\right).
\end{equation} 
Assuming that a Ramanujan--Petersson type bound for the coefficients of $L(\pi,s)$,  $\lambda_\pi(n)\ll n^{o(1)}$ holds, one obtains the trivial bound 
\begin{equation}\label{deltasum}
 \Delta(\lambda_\pi,X,a;q)\ll (qX)^{o(1)}\frac{X}q,
\end{equation}
and the question is to improve this bound for $q$ as large as possible relative to $X$.

The second term on the right hand side of \eqref{Deltadef} is easily evaluated in terms of the order of the pole of $L(\pi,s)$ as $s=1$. As for the first term
$$\mathop{\sum_{n\geq 1}}_{n\equiv a\bmod q}\lambda_\pi(n)V\left(\frac{n}{X}\right)=\sum_{n\geq 1}\lambda_\pi(n)\delta_{a\mods q}(n)V\left(\frac{n}{X}\right),$$
the Dirac function $\delta_{a\mods q}(n)$  is not a trace function in the above sense\footnote{Although the scaled function $q^{1/2}\delta_{a\mods q}(n)$ might reasonably be considered as such.} but an application of the functional equation for the character twists $L$-functions $L(\pi\times\chi,s)$ for $\chi$ varying over the Dirichlet characters of modulus $q$ transforms the left hand side of \eqref{Deltadef} into a sum essentially of the shape
\begin{equation}\label{deltatransform}
	\frac{X}{q^{\frac{d+1}2}}\sum_{n\geq 1}\lambda_{\ov\pi}(n)\Kl_d(an;q)\widecheck{V}\left(\frac{n}{q^d/X}\right).
\end{equation}
 Here $\widecheck{V}(x)$ denotes a rapidly decreasing function  which is a suitable integral transform of $V$ (depending on $d$ and the Gamma factors of $\pi$) and
$$\Kl_d(n; q)=\frac{1}{q^{\frac{d-1}{2}}}\sumsum_\stacksum{x_1,\cdots,x_d\in(\Zz/q\Zz)^\times}{x_1.\cdots.x_d=n}e\left(\frac{x_1+\cdots+x_d}{q}\right)$$
denotes the $d$-th hyper-Kloosterman sum. As is well known, Kloosterman sums are trace functions  (see \cite{GKM}) and they satisfy Deligne's bound
$$|\Kl_d(n; q)|\leq d.$$ 
Therefore, possibly subject to the Ramanujan--Petersson conjecture, one obtains that the dual sum can be bounded as 
$$\frac{X}{q^{\frac{d+1}2}}\sum_{n\geq 1}\lambda_{\ov\pi}(n)\Kl_d( a n;q)\widecheck{V}\left(\frac{n}{q^d/X}\right)\ll  (qX)^{o(1)}q^{\frac{d-1}{2}}.$$
This bound improves \eqref{deltasum} as long as
\begin{equation}
	\label{arithprogrrange}
	q\leq X^{\theta_d-\delta},\ \delta>0
\end{equation}
where \begin{equation}
	\label{thetaddef}
	\theta_d=\frac{2}{d+1}.
\end{equation}
We call the exponent $\theta_d$ the {\em standard level of distribution} of the sequence $(\lambda_\pi(n))_{n\geq 1}$ for individual (prime) moduli.
\begin{remark} For $d=2$ this reasoning is due to Selberg in the case of the divisor function while for higher values of $d$ it can be found in the work \cite{LRS} of Luo--Rudnick--Sarnak. 
\end{remark}

In the three cases described above we obtain
$$\theta_2=2/3,\ \theta_3=1/2,\ \theta_6=2/7.$$

A natural question is then whether one can enlarge the standard level of distribution $\theta_d$. Considering the limit case $$X=q^{1/\theta_d}=q^{\frac{d+1}2}$$ we see that it would amount to obtaining \eqref{nontrivialgeneric} for $K(n)=\Kl_d(an;q)$ in any range shorter than 
$$\widecheck{X}\asymp q^d/X=q^{\frac{d-1}2}.$$
This range, which we call the {\em a.p. range}, is shorter than the convexity range \eqref{convexrange} by a factor $q^{1/2}$.

\begin{remark}
The three results \cites{FKM1,KLMS,LMS} mentioned above, while improving over the convexity range, fall short of reaching the a.p. range.	
\end{remark}
 
Inspired by a recent work  \cite{huang} of Huang, we provide an example in degree $4$ for which the standard level of distribution $\theta_4=2/5$ can be improved:
\begin{theorem}\label{main-theorem}
Let $f$ be either a holomorphic Hecke eigencuspform or a Hecke--Maass cusp form of full level whose Hecke eigenvalues are denoted $(\lambda_f(n))_{n\geq 1}$ and let $\varpi_f\geq 0$ be such that the following bound holds
\begin{equation}
    \label{RPbound}
    |\lambda_f(n)|\leq n^{\varpi_f+o(1)},\ n\ra\infty
\end{equation}
and let
\begin{equation}
    \label{thetafdef}
    \theta_{f}=\frac{21}{52(1+4\varpi_f)}.
\end{equation}
Let $q$ be a prime and $a\geq 1$ be an integer such that $(a,q)=1$. 
For any $X\geq 1$  and $\eta>0$ satisfying 
\begin{equation}
    \label{qbound}
    q\leq X^{\theta_f-\eta},
\end{equation}
we have
	$$\mathop{\sum_{n\geq 1}}_{n\equiv a\bmod q}\lambda_f(n)^2V\left(\frac{n}X\right)-\frac{1}{\varphi(q)}\sum_\stacksum{n\geq 1}{(n,q)=1}\lambda_f(n)^2V\left(\frac{n}X\right)\ll_{f,V,\eta} (X/q)^{1-\delta}$$
	for some $\delta=\delta(\eta)>0$.
	In particular if $f$ is holomorphic $\varpi_f=0$ and 
	$$\theta_f=21/52=2/5+1/260.$$
\end{theorem}
\begin{remark}\label{rmksquare} The Ramanujan--Petersson conjecture for cusp forms predicts that
for any $f$, $$\varpi_f=0.$$
The best known bound, due to Kim and Sarnak \cite{KimSar}*{Proposition 2} is
$$\varpi_f\leq 7/64.$$
This is unfortunately not sufficient to insure that
$\theta_f>\theta_4=2/5$. For this we would need $7/64$ to be replaced by any exponent smaller than $1/416$. 

We also observe that the closely related Theorem \ref{main-theorem2} below is independent of any approximation toward the Ramanujan--Petersson conjecture.

In the work \cite{huang}, Huang improved the error term for the sharp-cut sum $$\sum_{n\leq X}\lambda_f(n)^2$$ from $O\big(X^{3/5}\big)$ to $O\big(X^{3/5-1/560+o(1)}\big)$, resolving a long standing problem going back to Rankin and Selberg. In that Archimedean case Huang was able to avoid the use of the pointwise bound \eqref{RPbound}.
\end{remark}

\subsection{Factorable arithmetic functions}

The proof of Theorem \ref{main-theorem} builds on the fact that the arithmetic function $n\mapsto\lambda_f(n)^2$ is a Dirichlet convolution.

Indeed, returning to the general problem \eqref{nontrivialgeneric}, if the arithmetic function $(\lambda_\pi(n))_{n\geq 1}$  is factorable, i.e., is the Dirichlet convolution of two (or more) arithmetic functions (themselves being coefficients of automorphic $L$-functions):
$$\lambda_\pi(n)=\lambda_{\pi_1}\star \lambda_{\pi_2}(n)=\sum_{lm=n}\lambda_{\pi_1}(l)\lambda_{\pi_2}(m)$$
then one can expect that the bilinear structure presented in the sum
$$S_V(K;X)=\sumsum_{l,m}\lambda_{\pi_1}(l)\lambda_{\pi_2}(m)K(lm;q)V(\frac{lm}{X}),$$ will help in improving the range for \eqref{nontrivialgeneric} as one would be able to apply different treatments to the variables $l$ and $m$: using Cauchy's inequality and Poisson summation formula and eventually the problem is reduced to bounds like \eqref{nontrivialgeneric} for different $X$'s and arithmetic functions of lower complexity. Here it is also a good occasion to recall the work \cite{Chen} of Chen, who introduced for the first time a {\em bilinear structure} to the reminder terms in the application of the linear sieve and non-trivially treated those reminder terms; Chen's idea in \cite{Chen} had inspired many far-reaching subsequent developments on sieve methods and their applications in the theory of prime numbers. 

The above principle is hardly  new and has already been exploited multiple times. A particularly striking example is the work of Friedlander--Iwaniec \cite{FrIwAnn} (see also \cites{HBActa,FKMMath} for subsequent improvements) who improved the standard level of distribution $\theta_3=1/2$,  of the ternary divisor function 
$$d_3(n)=1\star 1\star 1(n)=\sumsum_{klm=n}1.$$
In that case, the ultimate goal is to bound non-trivially the trilinear sum
$$\sumsum_{k,l,m}\Kl_3(aklm;q)V(\frac{klm}{X})$$
near the a.p. range $X=q$. This was achieved by different methods depending on the relative sizes of the variables $k,l,m$; the  critical case was
$$kl=q^{1/2-\eta+o(1)},\ m=q^{1/2-\eta+o(1)}$$
for $\eta>0$ and small enough; this roadblock was overcome by using a shifting technique dating back at least to Vinogradov and Korobov, supplemented by sharp bounds on rather complicated sums of products of Kloosterman sums due to  Birch and Bombieri building on Deligne's work on the Riemann Hypothesis over finite fields.

\begin{remark}
This work \cite{FrIwAnn} of Friedlander--Iwaniec was recently extended in \cites{KMS,KMS2} to Dirichlet convolutions of the form
$$1\star \lambda_f(n)=\sumsum_{lm=n}\lambda_f(m)$$
for $f$ a Hecke cuspform; again the critical range was 
$$l=q^{1/2-\eta+o(1)},\ m=q^{1/2-\eta+o(1)}.$$
Notice that $(d_3(n))_{n\geq 1}$ and $(1\star \lambda_f(n))_{n\geq 1}$ are the coefficients attached to $\GL_{3,\Qq}$ Eisenstein series representations (the isobaric sums $1\boxplus 1\boxplus 1$ and $1\boxplus \pi_f$, respectively) and by the same methods, one can improve the standard level of distribution $\theta_3=1/2$ of $(\lambda_\pi(n))_{n\geq 1}$ for $\pi$ any fixed $\GL_{3,\Qq}$ Eisenstein series representation. Doing the same for $\pi$ a cuspidal representation seems to be a serious challenge.
\end{remark}

To prove Theorem \ref{main-theorem} we use the fact that $n\mapsto \lambda_f(n)^2$ is (essentially) a Dirichlet convolution: we have the well known identity
$$L(f\times f,s):=\zeta(2s)\sum_{n\geq 1}\frac{\lambda_f(n)^2}{n^s}=\zeta(s)\sum_{m\geq 1}\frac{\lambda_{\sym_2f}(m)}{m^s}=:\zeta(s)L(\sym_2f,s)$$
where $\sym_2 f$ denotes the $\GL_{3,\Qq}$ symmetric square lift of (the automorphic representation attached to) $f$ and $L(\sym_2f,s)$ is the standard $L$-function associated to it. Therefore we obtain 
\begin{equation}
    \label{convident}
    \lambda_f(n)^2=\sum_{d^2k=n}\mu(d)\times \lambda_{1\boxplus{\sym_2f}}(k)=\sum_{d^2 lm=n}\mu(d)\lambda_{\sym_2f}(m).
\end{equation}
Theorem \ref{main-theorem} is then a simple consequence of
\begin{theorem}\label{main-theorem2}Let $\pi$ be a $\GL_{3,\Qq}$ cuspidal automorphic representation of conductor $1$ and $(\lambda_\pi(n))_{n\geq 1}$ be the coefficient of its standard $L$-function. Set  $$\lambda_{1\boxplus\pi}(n):=1\star\lambda_{\pi}(n):=\sum_{lm =n}\lambda_{\pi}(m).$$
Let $q$ be a prime and $a\geq 1$ be an integer such that $(a,q)=1$. For any $\eta>0$ and $X$ satisfying
\begin{equation}
    \label{qbound2}
    q\leq X^{21/52-\eta}=X^{2/5+\expo-\eta},
\end{equation}
we have 
$$\mathop{\sum_{n\geq 1}}_{n\equiv a\bmod q}\lambda_{1\boxplus\pi}(n)V\left(\frac{n}{X}\right)-\frac{1}{\varphi(q)}\sum_\stacksum{n\geq 1}{(n,q)=1}\lambda_{1\boxplus\pi}(n)V\left(\frac{n}{X}\right)\ll_{\pi,V,\eta} (X/q)^{1-\delta}$$
for some $\delta=\delta(\eta)>0$.
	
\end{theorem}

\begin{remark} Unlike \eqref{qbound} the exponent in  \eqref{qbound2} is independent of any approximation to the Ramanujan--Petersson conjecture for $\GL_{3}$ automorphic representations. 
\end{remark}

\section{Proof of Theorem \ref{main-theorem}}

In this section we deduce Theorem \ref{main-theorem} from Theorem \ref{main-theorem2}.

Set
$\pi=\sym_2f$ (which is cuspidal since $f$ has level $1$ and is therefore not CM).
Given $a,d\geq 1$ integers coprime with $q$, we denote by $a'=a'(a,d,q)$ any integer such that 
$$d^2a'\equiv a\mods q.$$
Using \eqref{Deltadef} and \eqref{convident}, we have
\begin{align*}
   \Delta(\lambda_f^2,X;a,q) &=\sum_{(d,q)=1}\mu(d)\Delta(\lambda_{1\boxplus\pi},X/d^2,a';q)\\&=\sum_{d\leq D}\mu(d)\Delta(\lambda_{1\boxplus\pi},X/d^2,a';q)+\sum_{d> D}\mu(d)\Delta(\lambda_{1\boxplus\pi},X/d^2,a';q),
\end{align*}
 where 
   $$D=X^{\eta'}$$ for some $\eta'>0$ small enough (we will choose $\eta'=\frac{2\varpi_f}{1+4\varpi_f}+\eta$ later) so that
   \begin{equation}\label{condition-on-q}
   q\leq (X/D^2)^{21/52-\eta/10}.
   \end{equation}
We bound the first sum $\sum_{d\leq D}(\cdots)$ by applying Theorem \ref{main-theorem2} to the inner term. Indeed, for $\delta<1/2$ we have
$$\sum_{d\leq D}\mu(d)\Delta(\lambda_{1\boxplus\pi},X/d^2,a';q)\ll \sum_{d\leq D}(\frac{X}{d^2q})^{1-\delta}\ll (X/q)^{1-\delta}.$$
For the second sum we apply the trivial bound
$$\Delta(\lambda_{1\boxplus\pi},X/d^2,a';q)\ll \frac{(X/d^2)\max_{n\ll X/d^2}|\lambda_{1\boxplus\pi}(n)|}{q}\ll X^{o(1)}\frac{(X/d^2)(X/d^2)^{2\varpi_f}}{q}$$
where we used $|\lambda_{\sym_2f}(n)|\ll n^{2\varpi_f+o(1)}$ with $\varpi_f$ satisfying \eqref{RPbound}. Since $$q\leq X^{\frac{21}{52(1+4\varpi_f)}-\eta}$$
we see that \eqref{condition-on-q} is satisfied
for 
$$D= X^{\frac{2\varpi_f}{1+4\varpi_f}+\eta}.$$ 
We obtain
$$\sum_{d>D}\mu(d)\Delta(\lambda_{1\boxplus\pi},X/d^2,a';q)\ll \frac{X^{1+o(1)}(X/D^2)^{2\varpi_f}}{Dq}\ll (X/q)^{1-\delta'}$$
for some $\delta'=\delta'(\eta)>0$.

Putting the two bounds together we conclude that
\begin{align*}
   \Delta(\lambda_f^2,X;a,q) \ll (X/q)^{1-\min(\delta,\delta')}.
\end{align*}

This completes the proof of Theorem \ref{main-theorem}.

\section{Proof of Theorem \ref{main-theorem2}}

More generally, let $\pi$ be a cuspidal automorphic representation of $\GL_{3,\Qq}$ of level $1$, we want to show that for $(a,q)=1$
$$\mathop{\sum_{n\geq 1}}_{n\equiv a\bmod q}\lambda_{1\boxplus\pi}(n)V(\frac{n}{X})=\frac{1}{\varphi(q)}\mathop{\sum_{n\geq 1}}_{(n,q)=1}\lambda_{1\boxplus\pi}(n)V(\frac{n}{X})+O_{\pi,V,\eta} \left((X/q)^{1-\delta}\right)$$
holds for $q<X^{2/5+\expo-\eta}$.

If we write $$K(n;q)=q^{1/2}\delta_{n\equiv a\bmod q},\ \widecheck{K}^4(n)=\mathrm{Kl}_4(an;q),$$
decompose $K(n;q)$ into a linear combination of Dirichlet characters $\chi\mods q$ and
 apply the functional equation for $L((1\boxplus\pi)\times \chi,s)$ in a way similar to \cite{LMS}*{Corollary 9.2}, we find that the left hand side above is equal to
\begin{multline}\label{Functional eq}
	q^{-1/2}\sum_{n\geq 1}\lambda_{1\boxplus\pi}(n)K(n;q)V(\frac{n}{X})=\frac{1}{\varphi(q)} \mathop{\sum_{n\geq 1}}_{(n,q)=1}\lambda_{1\boxplus\pi}(n)V(\frac{n}{X})\\
	+\frac{X}{q^{5/2}}\sum_{n\geq 1}\lambda_{1\boxplus\ov{\pi}}(n)\mathrm{Kl}_4(an;q)\widecheck{V}^4(\frac{n}{q^4/X})-\frac{X}{q^4\varphi(q)}\sum_{n\geq 1}\lambda_{1\boxplus\ov{\pi}}(n)\widecheck{V}^4(\frac{n}{q^4/X});
\end{multline}
see Lemma \ref{functional-equ} for a proof of this identity.
Here
$$\widecheck{V}^4(y)=\intc_{(3/2)}  \frac{L_\infty(1\boxplus\ov{\pi},s+\kappa)}{L_\infty(1\boxplus\pi,1-s+\kappa)}\widetilde V(1-s)y^{-s}ds$$
is a rapidly decreasing function of $y$ and  $\widetilde V(s)=\int_0^\infty V(y)y^{s}\frac{dy}{y}$ denotes the Mellin transform of $V$.

Plugging in the definition $$\lambda_{1\boxplus\ov{\pi}}(n)=\sum_{l m=n}\lambda_{\ov{\pi}}(m)$$ and applying the Rankin--Selberg bound $\sum_{m\leq X}|\lambda_{\ov{\pi}}(m)|^2\ll X$, it is easily seen that the last term on the right hand side of \eqref{Functional eq} contributes at most $O(q^{-1+o(1)})$. As for the second term, we have
\begin{equation*}
\begin{split}
\frac{X}{q^{5/2}}\sum_{n\geq 1}\lambda_{1\boxplus\ov{\pi}}(n)\mathrm{Kl}_4(an;q)\widecheck{V}^4(\frac{n}{q^4/X})=\frac{X}{q^{5/2}}\sum_{l,m\geq 1}\lambda_{\ov{\pi}}(m)\mathrm{Kl}_4(al m;q)\widecheck{V}^4(\frac{l m}{q^4/X}).
\end{split}
\end{equation*}
By introducing a dyadic partition of unity, we are reduced to considering sums of the form
\begin{equation}
    \label{MNsums}
    \frac{X}{q^{5/2}}\sum_{l\geq 1}\sum_{m\geq 1}\lambda_{\ov{\pi}}(m)\mathrm{Kl}_4(al m;q)V_1(\frac{l }{L})V_2(\frac{m}{M})
\end{equation}
for $O(\log^2 X)$ many real numbers $L,M\geq 1$ satisfying 
\begin{equation}
    \label{LMupperbound}
    LM\ll \frac{q^4}{X}.
\end{equation}
Since $|\Kl_4(alm;q)|\leq 4$ and $\sum_{m\leq M}|\lambda_{\ov{\pi}}(m)|^2\leq_{\ov{\pi}}M$ the trivial bound is
\begin{equation}
    \label{MNboundtrivial}
    \frac{X^{1+o(1)}}q(\frac{LM}{q^{3/2}}),
\end{equation}
which is good enough if $q\leq X^{2/5-\eta}$, and that henceforth one assumes that $q\geq X^{2/5-\eta}$. In particular, we may assume that
\begin{equation}
    \label{LMlower}
    LM\geq q^{3/2-\eta}
\end{equation}
for some fixed $\eta>0$ that can be chosen as small as necessary. In particular
$$L+M\geq q^{3/4-\eta}.$$

To obtain nontrivial cancellation for the sum \eqref{MNsums}, we split the argument into several cases. This strategy of the proof has  been somehow carried out in  \cite{KLMS}*{\S 10} with $\lambda_\pi(m)$ denoted $\lambda_\vphi(1,m)$. For completeness we include the details below.

 For any function $K(\cdot)$ on $\Zz/q\Zz$ we denote
\begin{equation}\label{FourierTransform}
\what K(n)=\frac{1}{q^{1/2}}\sum_{x\mods q}K(x)e(\frac{nx}q),\ e(\cdot)=\exp(2\pi i\cdot)
\end{equation}
its unitarily normalized Fourier transform; likewise, for $V$ a Schwartz function on $\Rr$, we denote its Fourier transform by
$$\what V(y)=\int_\Rr V(x)e(-xy)dx.$$
The treatment of these depends on the relative sizes of $L$ and $M$.

\subsection{The case $M\geq q^{4/3}$}
If $M$ is that long we may apply \cite{KLMS}*{Theorem 1.3} in two different ways depending on the size of $L$ compared to $q$.

If $L$ is small (say $L\leq q^{1/2}$) we apply
\cite{KLMS}*{Theorem 1.3} directly (note that $M\leq LM\leq q^{2}$) getting
\begin{multline}\label{KLMSbound1}
	\frac{X}{q^{5/2}}\sum_{l\geq 1}\Bigl(\sum_{m\geq 1}\lambda_{\ov{\pi}}(m)\mathrm{Kl}_4(al m;q)V_2(\frac{m}{M})\Bigr)V_1(\frac{l }{L})\\
	\ll_{\ov{\pi}} \frac{X^{1+o(1)}}{q^{5/2}}\sum_{l \geq 1}\|\widehat{\mathrm{Kl}_4}\|_{\infty} q^{2/9}M^{5/6}|V_1(\frac{l }{L})|\ll \frac{X^{1+o(1)}}{q}(\frac{L^3q^{37}}{X^{15}})^{1/18}.
\end{multline}
The last step follows from the identity
\begin{equation}\label{Kl4fourier}
    \widehat{\mathrm{Kl}_4}(al \bullet;q)(\tilde m)=
\delta_{(\tilde m,q)=1}
\mathrm{Kl}_3(-al\overline{\tilde m};q)+q^{-2}
\end{equation}
and $LM\ll q^4/X$.

In particular this bound is suitable as long as
$$q\leq X^{15/37-\eta}L^{-3/37}$$
for some fixed $\eta>0$.
\begin{remark}
In particular, since $L\geq 1$, this implies that $$q\leq X^{15/37}$$ (which is the limit of our method, fortunately $15/37>2/5$).

In view of this and \eqref{LMupperbound} we may assume that
$$LM\leq q^{4-37/15}=q^{23/15}$$
which implies (since we have assumed $M\geq q^{4/3}$) that $L\leq q^{1/5}$.
\end{remark}

\subsection{The case $L\geq q^{1/2}$}
 In that situation we can improve over the trivial bound by applying the Poisson summation formula in the $l$ variable: using \eqref{Kl4fourier}
 \begin{gather}\nonumber
     \frac{X}{q^{5/2}}\sum_{m\geq 1}\lambda_{\ov{\pi}}(m)\bigl(\sum_{l\geq 1}\Kl_4(al m;q)V_1(\frac{l}{L})\bigr)V_2(\frac{m}{M})\\\nonumber
     =\frac{X}{q^{5/2}}\sum_{m\geq 1}\lambda_{\ov{\pi}}(m)\left(\frac{L}{q^{1/2}}\sum_{\tilde{l}\in \Zz}\bigl(\delta_{(\tilde l,q)=1}\mathrm{Kl}_3(-am\overline{\tilde{l}};q)+q^{-2}\bigr)\widehat{V_1}(\frac{\tilde{l}}{q/L})\right)V_2(\frac{m}{M})
     \\
\ll\frac{X^{1+o(1)}}{q^{5/2}}Mq^{1/2}= \frac{X^{1+o(1)}}{q}\frac{M}q.\label{MNbound2}
 \end{gather}

This bound is good as long as
$$M\leq q^{1-\eta},\ \eta>0$$
which occurs as soon as
\begin{equation}
    \label{goodcase2}
    L\geq q^{8/15+\eta}.
\end{equation}

\subsection{The case $L\leq q^{1-\eta}$ ($\eta>0$)} We can also have some gain by applying Cauchy--Schwarz with $l$ inside the square, followed with Poisson in the $m$-variable. Indeed, we have
\begin{multline*}
	\frac{X}{q^{5/2}}\sum_{m\geq 1}\lambda_{\ov{\pi}}(m)\bigl(\sum_{l\geq 1}\mathrm{Kl}_4(al m;q)V_1(\frac{l}{L})\bigr)V_2(\frac{m}{M})\ll \\ \frac{X}{q^{5/2}}\left(\sum_{m\geq 1}|\lambda_{\ov{\pi}}(m)|^2V_2(\frac{m}{M})\right)^{1/2}\left(\sum_{m\geq 1}|\sum_{l\geq 1}\mathrm{Kl}_4(al m;q)V_1(\frac{l}{L})|^2V_2(\frac{m}{M})\right)^{1/2}\end{multline*}
\begin{multline}\\\ \ll\frac{XM^{1/2}}{q^{5/2}}\left(\sum_{l_1,l_2\geq 1}V_1(\frac{l_1}{L})\overline{V_1(\frac{l_2}{L})}\sum_{m\geq 1}\mathrm{Kl}_4(al_1 m;q)\overline{\mathrm{Kl}_4(al_2 m;q)}V_2(\frac{m}{M})\right)^{1/2}.\label{c-sum}
\end{multline}

We consider two subcases.
\subsubsection{If $l_1= l_2:=l$} The sum inside the parentheses above can be simply bounded by
$$\sum_{l\geq 1}|V_1(\frac{l}{L})|^2\sum_{m\geq 1}|\mathrm{Kl}_4(al m;q)|^2V_2(\frac{m}{M})\ll LM.$$
\subsubsection{If $l_1\not= l_2$}

We apply Poisson summation in the $m$-variable, getting
$$\sum_{m\geq 1}\mathrm{Kl}_4(al_1 m;q)\overline{\mathrm{Kl}_4(al_2 m;q)}V_2(\frac{m}{M})=\frac{M}{q^{1/2}}\sum_{\tilde{m}\in \Zz}\mathcal{C}_{a}(\tilde{m},l_1,l_2;q)\widehat{V_2}(\frac{\tilde{m}}{q/M})$$
where
$$\mathcal{C}_{a}(\tilde{m},l_1,l_2;q):=\frac{1}{q^{1/2}}\sum_{x\in(\Zz/q\Zz)^\times}\mathrm{Kl}_4(al_1 x;q)\overline{\mathrm{Kl}_4(al_2 x;q)}e(\frac{x \tilde{m}}{q}).$$
By \cite{FKM-sumproduct}*{Corollary 3.2}, we have the following bound 
$$\mathcal{C}_{a}(\tilde{m},l_1,l_2;q)\ll q^{1/2}\delta_{\substack{\tilde{m}\equiv 0\mods{q}\\ l_1\equiv l_2\mods{q}}}+1.$$
Since we have assumed that $l_1,l_2<q$ we also have $l_1\not\equiv l_2\mods q$ and
the $m$-sum is bounded by
$$
    \frac{M}{q^{1/2}}\sum_{\tilde{m}}\bigl|\mathcal{C}_{a}(\tilde{m},l_1,l_2;q)\widehat{V_2}(\frac{\tilde{m}}{q/M})\bigr|
\ll q^{1/2+o(1)}
$$
and the original sum can be bounded as follows
\begin{align}\nonumber
    	\frac{X}{q^{5/2}}\sum_{l\geq 1}\sum_{m\geq 1}\lambda_{\ov{\pi}}(m)\mathrm{Kl}_4(al m;q)V_1(\frac{l}{L})V_2(\frac{m}{M})&\ll  \frac{X^{1+o(1)}M^{1/2}}{q^{5/2}}\bigl(LM+L^2q^{1/2}\bigr)^{1/2}\\
&\ll
\frac{X^{1+o(1)}}q\bigl(\frac{LM^2}{q^3}+\frac{L^2M}{q^{5/2}}\bigr)^{1/2}\nonumber \\
&\ll
\frac{X^{1+o(1)}}q\bigl(\frac{1}{L}\frac{q^5}{X^2}+L\frac{q^{3/2}}{X}\bigr)^{1/2}.\label{MNbound3}
\end{align}
In view of \eqref{goodcase2} we will apply this bound only when 
$$L\leq q^{8/15+\eta}$$ for $\eta>0$ small enough (in particular so that $q^{8/15+\eta}\leq q^{ 1-\eta}$). Assuming this,  the second term in the parentheses on the right hand side of  \eqref{MNbound3} satisfies
$$L\frac{q^{3/2}}{X}\leq X^{-13/74+\eta}.$$ 
Therefore, under these conditions \eqref{MNbound3} is good as soon as
\begin{equation}
    \label{goodcase3}
  q\leq X^{2/5-\eta}L^{1/5}.
\end{equation}
\subsection{Putting it all together}
Let $$L_0=X^{1/52}$$ be the solution of the equation $$X^{15/37}L_0^{-3/37}=X^{2/5}L_0^{1/5}=X^{21/52}=X^{2/5+\expo}.$$

We need to show that for any  small enough  $\eta>0$ and any prime $q$ satisfying
$$X^{2/5-\eta}\leq q\leq X^{21/52-\eta},$$ one has
$$\Delta(\lambda_{1\boxplus\pi},X,a;q)\ll (X/q)^{1-\delta}$$
for some $\delta=\delta(\eta)>0$. 

It is sufficient to show that this bound holds for any of the sums \eqref{MNsums} for $L,M$ satisfying
$$1\leq LM\leq q^4/X.$$
\begin{itemize}
   \item[-] If $M\leq  q^{4/3}$ and  $L\leq L_0$ then $LM\leq q^{3/2-\delta}$ for some $\delta=\delta(\eta)>0$ and we just use the trivial bound \eqref{MNboundtrivial}.
   \item[-] If  $M\leq  q^{4/3}$ and  $L_0\leq L\leq q^{8/15+\eta}$ we use \eqref{MNbound3}.
   \item[-] If  $M\leq  q^{4/3}$ and  $L\geq q^{8/15+\eta}$ we use \eqref{MNbound2}.
    \item[-] If $M\geq q^{4/3}$ and $L\leq L_0$, we use \eqref{KLMSbound1}.
     \item[-] If $M\geq q^{4/3}$ and $L\geq L_0$, then $L\leq q^{1/5}$ and we use again \eqref{MNbound3}.
\end{itemize}

\begin{remark} The above proof can be carried out for more general convolutions
$1\star\lambda_\pi$
for $\lambda_\pi$ of degree $d\geq 2$ such that \eqref{nontrivialgeneric} can be obtained near and below the convexity range $q^{d/2}$: let $X=q^{\frac{d+2}2}$ be the a.p. range for $1\star\lambda_\pi$. By duality and a dyadic decomposition, one has to bound non-trivially  bilinear sums of the shape
$$\sum_{l\sim L}\sum_{m\sim M}\lambda_{\ov\pi}(m)\Kl_{d+1}(lm;q)$$
for $L,M\geq 1$ such that $LM\leq q^{d/2}$. Let us assume for simplicity that $LM=q^{d/2}$.
\begin{enumerate}
	\item If $L\leq q^{\eta}$ for $\eta>0$ small enough then the $m$-sum 
$$\sum_{m\sim M}\lambda_{\ov\pi}(m)\Kl_{d+1}(lm;q)$$
can be bounded non-trivially since $M\geq q^{d/2-\eta}$ is near the convexity range for $\lambda_{\ov\pi}$. 

\item If $L\geq q^{d/2-\eta}$ the linear sum
$$\sum_{l\sim L}\Kl_{d+1}(lm;q)$$
can be bounded non-trivially after applying the Poisson summation formula (as $d/2-\eta>1/2$ because $d\geq 2$). 

\item In the remaining range
we have
$$\min(L;M)\geq q^{\eta}, \max(L,M)\geq q^{\frac{d}4}.$$
\begin{itemize}
    \item  If $d\geq 3$ we can then apply Cauchy--Schwarz to this bilinear sum with the shorter variable inside the square and the longer variable being smoothed;
 since the longer variable is well above the P\'olya--Vinogradov range  $q^{1/2}$ (since $q^{d/4}\geq q^{3/4}$) there will be some saving from sums of products of Kloosterman sums and we are done.
 \item In the case $d=2$, then in the last situation we obtain only 
 $\max(L,M)\geq q^{1/2}$ which is the P\'olya--Vinogradov range and the above approach does not provide any saving if $L=M=q^{1/2}$. However using a technique going back to Vinogradov--Korobov  developed in this context by Friedlander--Iwaniec in \cite{FrIwAnn} one can still bound non-trivially the bilinear form of Kloosterman sums and eventually get some saving (see \cites{KMS,KMS2}).
\end{itemize}	
\end{enumerate} 
For instance the above method can be applied to the convolution of degree $1+6=7$   $$1\star\lambda_{\pi}$$ 
to pass above the standard distribution exponent $\theta_7=2/(7+1)=1/4$
for $1\star\lambda_{\pi}$ with
$\lambda_{\pi}(n)$
given by \eqref{RS23def}. The main point is that  the Kloosterman sheaf $\KL_7$ is a ``good" sheaf
in the sense of \cite{LMS}*{\S 1} so that \eqref{nontrivialgeneric} holds for $K(\bullet)=\Kl_7(l\bullet;q)$ at or slightly below the convexity range $q^{6/2}=q^3$ (see \cite{LMS}*{Theorem 1.1}).
\end{remark}

\begin{lemma}\label{functional-equ}
Proof of the functional equation \eqref{Functional eq}.
\end{lemma}
\proof 
The proof is exactly the same as that of \cite{LMS}*{Corollary 9.2}.
We denote  $\kappa=0$ if $\chi(-1)=1$ and $\kappa=1$ if $\chi(-1)=-1$ and denote
 $$\eps_\chi=q^{-1/2}\sum_{x\in\Fqt}\chi(x)e(\frac{x}q)$$
the normalized Gauss sum. 

Recall that  $$\lambda_{1\boxplus\pi}(n):=1\star\lambda_{\pi}(n):=\sum_{lm =n}\lambda_{\pi}(m).$$

The $L$-function $$L((1\boxplus\pi)\times \chi,s)=\sum_{n\geq 1}\frac{\lambda_{1\boxplus\pi}(n)\chi(n)}{n^s}=L(\chi,s)L(\pi\times \chi,s)$$ has analytic continuation to $\Cc$ and satisfies a functional equation of the form:
\begin{equation}\label{FE-gl4}
\Lambda((1\boxplus\pi)\times \chi,s)=\eps_\chi^4\Lambda((1\boxplus\ov{\pi})\times \overline{\chi},1-s)\end{equation}
where
$$\Lambda((1\boxplus\pi)\times \chi,s)=q^{2s}L_\infty(1\boxplus\pi,s+\kappa)L((1\boxplus\pi)\times \chi,s)$$
is the completed $L$-function and 
 $$L_\infty(1\boxplus\pi,s)=\prod_{i=1}^4\Gamma_\Rr(s-\mu_{i}),\ \Gamma_\Rr(s)=\pi^{-s/2}\Gamma(s/2)$$
with
$$\{\mu_{i},\ i=1,2,3,4\}=\{0,\nu_2-\nu_1,2\nu_1+\nu_2-1,1-\nu_1-2\nu_2\}$$
denoting the local Archimedean factor of $1\boxplus\pi$.

We have 
$$V(x)=\intc \widetilde V(s)x^{-s}ds$$
(the integration is along the vertical line $\Re s=1+1/14$) so that
\begin{equation*}
\begin{split}
&\mathop{\sum_{n\geq 1}}_{n\equiv a\bmod q}\lambda_{1\boxplus\pi}(n)V(\frac{n}{X})=\frac{1}{\varphi(q)}\sum_{\chi\mods q}\overline{\chi}(a)\sum_{n\geq 1}\lambda_{1\boxplus\pi}(n)\chi(n)V(\frac{n}{X})\\
=&\frac{1}{\varphi(q)} \mathop{\sum_{n\geq 1}}_{(n,q)=1}\lambda_{1\boxplus\pi}(n)V(\frac{n}{X})+ \frac{1}{\varphi(q)}\sum_\stacksum{\chi\mods q}{\chi\not=\chi_0}\overline{\chi}(a)\intc\widetilde{V}(s)\frac{\Lambda((1\boxplus\pi)\times \chi,s)}{L_\infty(1\boxplus\pi,s+\kappa)}(\frac{X}{q^2})^sds.
\end{split}\end{equation*}
By applying the functional equation \eqref{FE-gl4}, the second term above equals
\begin{equation*}
\begin{split}
\frac{q^2}{\varphi(q)}\sum_\stacksum{\chi\mods q}{\chi\not=\chi_0}\overline{\chi}(a)\eps_\chi^4\,\intc\widetilde{V}(s)L((1\boxplus\ov{\pi})\times \overline{\chi},1-s)\frac{L_\infty(1\boxplus\ov{\pi},1-s+\kappa)}{L_\infty(1\boxplus\pi,s+\kappa)}(\frac{X}{q^4})^sds.
\end{split}\end{equation*}

In the integral we make the change of variable $s\leftrightarrow 1-s$ getting
$$\frac{X}{q^2\varphi(q)}\sum_\stacksum{\chi\mods q}{\chi\not=\chi_0}\overline{\chi}(a)\eps_\chi^4\intc_{(-1/14)}  L((1\boxplus\ov{\pi})\times \overline{\chi},s)\frac{L_\infty(1\boxplus\ov{\pi},s+\kappa)}{L_\infty(1\boxplus\pi,1-s+\kappa)}\widetilde V(1-s)(\frac{X}{q^4})^{-s}ds$$
and shifting the contour back to $\Re s=3/2$ without hitting any poles we obtain the sum
\begin{equation*}
\begin{split}
&\frac{X}{q^2\varphi(q)}\sum_{n\geq 1}\lambda_{1\boxplus\ov{\pi}}(n)\bigl(\sum_\stacksum{\chi\mods q}{\chi\not=\chi_0}\overline{\chi}(an)\eps_\chi^4\bigr)\intc_{(3/2)}  \frac{L_\infty(1\boxplus\ov{\pi},s+\kappa)}{L_\infty(1\boxplus\pi,1-s+\kappa)}\widetilde V(1-s)(\frac{nX}{q^4})^{-s}ds\\
=&\frac{X}{q^2\varphi(q)}\sum_{n\geq 1}\lambda_{1\boxplus\ov{\pi}}(n)\bigl(\varphi(q)q^{-1/2}\mathrm{Kl}_4(an;q)-q^{-2}\bigr)\widecheck{V}^4(\frac{n}{q^4/X}).
\end{split}\end{equation*}
This completes the proof of \eqref{Functional eq}.

 \qed
 
 \subsection*{Acknowledgements} We would like to thank Bingrong Huang for sharing his work \cite{huang} which gave us the motivation for writing up this note. We thank the referees for their helpful suggestions.
 
 Chen's $2N=p+P_2$ Theorem is one of the greatest classics of the Analytic Number Theory literature and, as students, it gave us the impetus  to pursue in this area. It is therefore a great pleasure to dedicate this work to the memory of Jingrun Chen on the occasion of the 50th anniversary of his landmark theorem.

\end{document}